\documentclass[12pt]{article}
\usepackage{amssymb}

\newtheorem{definition}{Definition}[section]
\newtheorem{theorem}[definition]{Theorem}
\newtheorem{lemma}[definition]{Lemma}
\newtheorem{corollary}[definition]{Corollary}
\newtheorem{remark}[definition]{Remark}
\newtheorem{example}[definition]{Example}

\newtheorem{problem}[definition]{Problem}
\newtheorem{note}[definition]{Note}

\def\R{\mathbb R}

\def\K{\mathbb K}
\def\D{\mathbb D}
\def\Z{\mathbb Z}

\newcommand{\beast}{\begin{eqnarray*}}
\newcommand{\eeast}{\end{eqnarray*}}

\begin{document}
\newenvironment{proof}{\noindent{\it Proof\/}:}{\par\noindent $\Box$\par}

\title{ \bf Two relations that generalize \\
the $q$-Serre relations and the \\
Dolan-Grady relations\footnote{
{\bf Keywords}. $q$-Racah polynomial,  Leonard pair,
Tridiagonal pair,
 Askey scheme, Askey-Wilson algebra, Subconstituent algebra,
Terwilliger algebra, Askey-Wilson polynomials, Serre relations,
Dolan-Grady relations,
quadratic algebra. \hfil\break
\noindent {\bf 2000 Mathematics Subject Classification}. 
05E30, 05E35, 33C45, 33D45. 
}}
\author{Paul Terwilliger  
}
\date{}
\maketitle
\begin{abstract} We define an algebra on two generators
which we call the {\it Tridiagonal algebra}, and we   
consider its irreducible modules. The algebra is defined
as follows. 
Let $\K$ denote a field, and let
  $\beta, \gamma, \gamma^*, \varrho, \varrho^*$ denote a sequence
  of scalars taken from
$\K$. The corresponding Tridiagonal algebra $T$ is  
 the associative $\K$-algebra with 1 generated by two 
symbols 
 $A$, $A^*$
subject to the  relations  
\beast
\lbrack A,A^2A^*-\beta AA^*A + 
A^*A^2 -\gamma (AA^*+A^*A)-\varrho A^*\rbrack &=& 0,
\\
\lbrack A^*,A^{*2}A-\beta A^*AA^* + AA^{*2} -\gamma^* (A^*A+AA^*)-
\varrho^* A\rbrack &=& 0,  
\eeast
where $\lbrack r,s\rbrack $ means $rs-sr$. We call these  
relations the {\it Tridiagonal relations}.
For 
$\beta = q+q^{-1}$, $\gamma = \gamma^*=0$, $\varrho=\varrho^*=0$,
the Tridiagonal relations are the $q$-Serre relations
\beast 
 A^3A^* - \lbrack 3\rbrack_q  A^2A^*A +  
 \lbrack 3\rbrack_q AA^*A^2  - A^*A^3   &=& 0,
\\
 A^{*3}A -  \lbrack 3\rbrack_q  A^{*2}AA^* +  
 \lbrack 3\rbrack_q A^*AA^{*2}  - AA^{*3} &=& 0,
\eeast
where $\lbrack 3 \rbrack_q = q + q^{-1}+1$.
For 
$\beta = 2$, $\gamma = \gamma^*=0$, $\varrho=b^2$, $\varrho^*=b^{*2}$,
the Tridiagonal relations are the Dolan-Grady relations
\beast
\lbrack A, \lbrack A, \lbrack A, A^* \rbrack \rbrack \rbrack &=&  b^2\lbrack
A, A^* \rbrack,
\\
\lbrack A^*, \lbrack A^*, \lbrack A^*, A \rbrack \rbrack \rbrack &=& 
b^{*2}\lbrack
A^*, A \rbrack.
\eeast
In the first part of this paper, we survey what is known
about irreducible
finite dimensional $T$-modules. We focus on how these modules
are related to 
 the Leonard pairs recently introduced by the present
author, and the more general Tridiagonal pairs recently introduced
by Ito, Tanabe, and the present author.  
In the second part of the paper, we construct an infinite dimensional
irreducible $T$-module based on the Askey-Wilson polynomials. This
module is on the vector space $\K\lbrack x \rbrack $
consisting of all polynomials in an indeterminant $x$ that have 
coefficients in $\K$. Denoting
by $A$ the linear transformation on 
 $\K\lbrack x \rbrack $ which is multiplication by $x$, and
denoting by $A^*$ an 
  Askey-Wilson
second order $q$-difference operator for $x$,
we show $A$ and $A^*$ satisfy a pair of Tridiagonal relations. Using
this we give   
 $\K\lbrack x \rbrack $ the structure of an irreducible $T$-module.
 The Askey-Wilson polynomials form 
 a basis for this module, and these basis elements
 are eigenvectors for $A^*$.
\end{abstract}

%

\section{Leonard pairs}
\medskip
\noindent Throughout this paper, $\K$ will denote an arbitrary
field.

\medskip
\noindent 
We begin by recalling the notion of a Leonard pair.

\begin{definition} \cite{LS99}
\label{def:lprecall}
Let 
 $V$ denote a  
vector space over $\K$ with finite positive dimension.
By a {\it Leonard pair} on $V$,
we mean an ordered pair $A, A^*$, where
$A:V\rightarrow V$ and $A^*:V\rightarrow V$ are linear transformations 
that 
 satisfy both (i), (ii) below. 
\begin{enumerate}
\item There exists a basis for $V$ with respect to which
the matrix representing $A$ is diagonal and the matrix
representing $A^*$ is irreducible tridiagonal.
\item There exists a basis for $V$ with respect to which
the matrix representing $A^*$ is diagonal and the matrix
representing $A$ is irreducible tridiagonal.

\end{enumerate}
(A tridiagonal matrix is said to be irreducible
whenever all entries immediately above and below the main
diagonal are nonzero).

\end{definition}

\begin{note} According to a common notational convention, for
a linear transformation $A$ the conjugate-transpose of $A$ is denoted
$A^*$. We emphasize we are not using this convention. In a Leonard
pair $A,A^*$, the linear transformations $A$ and $A^*$
are arbitrary subject
to (i),  (ii) above.
\end{note}
\noindent
Here is an example of a Leonard pair.
Set 
$V={\K}^4$ (column vectors), set 
\beast
A = 
\left(
\begin{array}{ c c c c }
0 & 3  &  0    & 0  \\
1 & 0  &  2   &  0    \\
0  & 2  & 0   & 1 \\
0  & 0  & 3  & 0 \\
\end{array}
\right), \qquad  
A^* = 
\left(
\begin{array}{ c c c c }
3 & 0  &  0    & 0  \\
0 & 1  &  0   &  0    \\
0  & 0  & -1   & 0 \\
0  & 0  & 0  & -3 \\
\end{array}
\right),
\eeast
and view $A$ and $A^*$  as linear transformations from $V$ to $V$.
We assume 
the characteristic of $\K$ is not 2 or 3, to ensure
$A$ is irreducible.
Then $A, A^*$ is a Leonard
pair on $V$. 
Indeed, 
condition (ii) in Definition
\ref{def:lprecall}
is satisfied by the basis for $V$
consisting of the columns of the 4 by 4 identity matrix.
To verify condition (i), we display an invertible  matrix  
$P$ such that 
$P^{-1}AP$ is 
diagonal and 
$P^{-1}A^*P$ is
irreducible tridiagonal.
Set 
\beast
P = 
\left(
\begin{array}{ c c c c}
1 & 3  &  3    &  1 \\
1 & 1  &  -1    &  -1\\
1  & -1  & -1  & 1  \\
1  & -3  & 3  & -1 \\
\end{array}
\right).
\eeast
 By matrix multiplication $P^2=8I$, where $I$ denotes the identity,   
so $P^{-1}$ exists. Also by matrix multiplication,    
\begin{equation}
AP = PA^*.
\label{eq:apeq}
\end{equation}
Apparently
$P^{-1}AP$ equals $A^*$ and is therefore diagonal.
By (\ref{eq:apeq}), and since $P^{-1}$ is
a scalar multiple of $P$, we find
$P^{-1}A^*P$ equals $A$ and is therefore irreducible tridiagonal.  Now 
condition (i) of  Definition 
\ref{def:lprecall}
is satisfied
by the basis for $V$ consisting of the columns of $P$. 

\medskip
\noindent The above example is a member of the following infinite
family of Leonard pairs.
For any nonnegative integer $d$, 
the pair
\begin{equation}
A = 
\left(
\begin{array}{ c c c c c c}
0 & d  &      &      &   &{\bf 0} \\
1 & 0  &  d-1   &      &   &  \\
  & 2  &  \cdot    & \cdot  &   & \\
  &   & \cdot     & \cdot  & \cdot   & \\
  &   &           &  \cdot & \cdot & 1 \\
{\bf 0} &   &   &   & d & 0  
\end{array}
\right),
\qquad A^*= \hbox{diag}(d, d-2, d-4, \ldots, -d)
\label{eq:fam1}
\end{equation}
is a Leonard pair on the vector space $\K^{d+1} $,
provided the 
 characteristic of $\K$ is zero or  greater than $d$.
This can be  proved by modifying the 
 proof for $d=3$ given above. One shows  
$P^2=2^dI$  and $AP= PA^*$, where 
$P$ denotes the matrix with $ij$ entry
\begin{equation}
P_{ij} =  
\Biggl({{ d }\atop {j}}\Biggr) {{}_2}F_1\Biggl({{-i, -j}\atop {-d}};2\Biggr)
\qquad \qquad (0 \leq i,j\leq d).
\label{eq:ex1}
\end{equation}
We follow the standard notation for
hypergeometric series given in \cite{gasperrahmanbk}.

\medskip
\noindent  In \cite{LS99} we obtain a complete classification
of Leonard pairs over an arbitrary field $\K$. We will not
need the details here, 
but we wish to discuss one aspect. There is a connection between
Leonard pairs and certain orthogonal polynomials contained
in the Askey scheme \cite{KoeSwa}.
Observe the ${{}_2F_1}$ that appears in
(\ref{eq:ex1}) 
is a Krawtchouk polynomial \cite{KoeSwa}.
There exist 
families of Leonard pairs similar to the one above in 
which 
 the Krawtchouk polynomial is replaced by one of the following.

\medskip
\centerline{
\begin{tabular}[t]{c|c}
        type & polynomial \\ \hline 
 ${{}_4F_3}$ & Racah \\ 
 ${{}_3F_2}$ & Hahn, dual Hahn \\ 
 ${{}_2F_1}$ & Krawtchouk \\ 
 ${{}_4\phi_3}$ & $q$-Racah \\ 
 ${{}_3\phi_2}$ & $q$-Hahn, dual $q$-Hahn \\ 
 ${{}_2\phi_1}$ & $q$-Krawtchouk (classical, affine, quantum, dual) 
\end{tabular}}

\medskip
\noindent
The above polynomials are defined in Koekoek and Swarttouw \cite{KoeSwa},
and the  connection to Leonard pairs  is given in 
\cite[ch. 15]{LS99} and 
\cite[p. 260]{BanIto}.
The above examples exhaust essentially all Leonard pairs in
the following sense. As we will see below in the section on
eigenvalues, associated with any  Leonard pair is a 
certain scalar $q$.
The above examples exhaust all Leonard pairs for which
 $q\not=-1$. 

\medskip
\noindent  
 There is a theorem due to  Leonard
 \cite{Leodual} and Bannai and Ito \cite[p. 260]{BanIto} 
 that gives a characterization
 of the 
polynomials from the
 above table for $\K=\R$.
 This  result has come to be known as Leonard's theorem.  
 Our above  mentioned 
classification of Leonard pairs
  amounts to a ``linear algebraic version'' of
 Leonard's theorem. 

\medskip
\noindent  In this paper we are concerned with
the following feature of Leonard pairs.
For the  Leonard pair
in (\ref{eq:fam1})  one can show 
\begin{eqnarray}
A^2A^* - 2AA^*A + A^*A^2 &=& 4A^*,  
\label{eq:tdpre1}
\\
A^{*2}A - 2A^*AA^* + AA^{*2} &=& 4A.
\label{eq:tdpre2}
\end{eqnarray}
This phenomenon
is not a coincidence. As we will see, 
every Leonard pair satisfies two
polynomial relations reminiscent of
(\ref{eq:tdpre1}),
(\ref{eq:tdpre2}).
Before pursuing this further, we 
consider 
 a certain generalization  of a Leonard pair.

\section{Tridiagonal pairs}

\medskip In this section we shift our attention
from a Leonard pair to a more general object called
a {\it Tridiagonal pair}. To define this we use the 
following notation.
Let $V$ denote
a vector space over $\K$ with finite positive dimension.
 Let $A:V\rightarrow V$ denote a linear transformation.
 By an {\it eigenspace} of $A$ we mean 
a nonzero  subspace of $V$ of the form
\beast
\lbrace v \in V \;\vert \;Av = \theta v\rbrace,
\eeast
where $\theta \in \K$.
The scalar $\theta$ is the associated eigenvalue.
We say $A$ is {\it diagonalizable} on $V$ whenever
$V$ is spanned by the eigenspaces of $A$.

\begin{definition} \cite{TD00}
\label{def:tdrecall}
Let $V$ denote
a vector space over $\K$ with finite positive dimension.
By a {\it Tridiagonal pair} (or {\it TD pair}) on $V$
we mean an ordered pair $A, A^*$, where
$A:V\rightarrow V$ and $A^*:V\rightarrow V$ are linear transformations 
satisfying (i)--(iv) below.
\begin{enumerate}
\item $A$ and $A^*$ are both diagonalizable on $V$.
\item There exists an ordering $V_0, V_1,\ldots, V_d$ of the  
eigenspaces of $A$ such that 
\begin{equation}
A^* V_i \subseteq V_{i-1} + V_i+ V_{i+1} \qquad \qquad (0 \leq i \leq d),
\label{eq:tdrecall1}
\end{equation}
where $V_{-1} = 0$, $V_{d+1}= 0$.
\item There exists an ordering $V^*_0, V^*_1,\ldots, V^*_\delta$ of
the  
eigenspaces of $A^*$ such that 
\begin{equation}
A V^*_i \subseteq V^*_{i-1} + V^*_i+ V^*_{i+1} \qquad \qquad (0 \leq i \leq \delta),
\label{eq:tdrecall2}
\end{equation}
where $V^*_{-1} = 0$, $V^*_{\delta+1}= 0$.
\item There is no subspace $W$ of $V$ such  that  both $AW\subseteq W$,
$A^*W\subseteq W$, other than $W=0$ and $W=V$.
\end{enumerate}
We say the above TD pair is over $\K$.
\end{definition}

\medskip
\noindent We now show any Leonard pair is a Tridiagonal pair.

\begin{lemma}
\label{lem:lpvstd}
Let $V$ denote a vector space over $\K$ with finite positive
dimension. Let 
$A:V\rightarrow V$ and $A^*:V\rightarrow V$ denote linear transformations.
Then the following are equivalent.
\begin{enumerate}
\item  $A,A^*$ is a Leonard pair on $V$.
\item $A,A^*$ is a TD pair on $V$, and for each of $A$, $A^*$
all eigenspaces have dimension 1.
\end{enumerate}
\end{lemma}

\begin{proof} 
$(i)\rightarrow (ii)$
We show $A$ and $A^*$ satisfy the conditions (i)--(iv) of 
Definition
\ref{def:tdrecall}. Condition (i) holds, since  
it is immediate from 
Definition \ref{def:lprecall} that $A$ and $A^*$ are diagonalizable
on $V$.
 We now show $A$ and $A^*$ satisfy Definition 
\ref{def:tdrecall}(ii). 
Let $v_0, v_1, \ldots, v_d$ denote the basis
for $V$ referred to in Definition
\ref{def:lprecall}(i).
 For $0 \leq i \leq d$, let $V_i$
denote the subspace of $V$ spanned by $v_i$. By
\cite[Lem. 1.3]{LS99} 
the eigenvalues
of $A$ associated with $v_0, v_1, \ldots, v_d$ are mutually
distinct, so $V_0, V_1, \ldots, V_d$ is an ordering of the
eigenspaces  of $A$. This ordering satisfies
(\ref{eq:tdrecall1}), since 
the matrix representing $A^*$ with respect to 
 $v_0, v_1, \ldots, v_d$ is tridiagonal.
We have now shown $A, A^*$ satisfy 
 Definition 
\ref{def:tdrecall}(ii). Interchanging the roles of $A$ and $A^*$, we
find 
$A, A^*$ satisfy 
 Definition 
\ref{def:tdrecall}(iii). It is immediate from \cite[Lem. 3.3]{LS99} 
that $A, A^*$ satisfy Definition 
\ref{def:tdrecall}(iv).  We have now shown $A,A^*$
satisfy the conditions (i)--(iv) of 
Definition
\ref{def:tdrecall},
so $A, A^*$  
is a TD pair
on $V$. From our above  comments the eigenspaces of $A$ all have
dimension 1. Similarly  
 the eigenspaces of $A^*$ all have
dimension 1.   

\noindent 
$(ii)\rightarrow (i)$ 
We show $A$ and $A^*$ satisfy 
conditions (i), (ii) of Definition \ref{def:lprecall}. Concerning 
condition (i), 
let $V_0, V_1, \ldots, V_d$ denote the ordering
of the eigenspaces of $A$ referred to in Definition
\ref{def:tdrecall}(ii).  By assumption each of
these eigenspaces has  dimension 1.
For $0 \leq i \leq d$,
let $v_i$ denote a nonzero element in $V_i$,
 and observe $v_0, v_1, \ldots, v_d$ is a 
basis for $V$. From the construction
the matrix representing $A$ with respect to 
this basis is diagonal. Using (\ref{eq:tdrecall1}),
we find the matrix representing $A^*$ with respect to 
this basis is tridiagonal. This tridiagonal matrix
is irreducible in view of 
Definition \ref{def:tdrecall}(iv). We have now shown
 $A$ and $A^*$ satisfy 
Definition \ref{def:lprecall}(i). Interchanging the roles of
$A$ and $A^*$, we find
 $A$ and $A^*$ satisfy 
Definition \ref{def:lprecall}(ii). 
It follows $A,A^*$ is a Leonard pair on $V$.
\end{proof}

%
\noindent Given the connection between Leonard pairs and the
polynomials of the Askey scheme, we find it worthwhile to investigate
TD pairs. We want to know if TD pairs correspond to multi-variable
generalizations of the above-mentioned polynomials, and if not,
 what they do correspond to. TD pairs have a lot
of structure and we believe a complete classification is possible.
At present we have nothing of the sort so we pose the following 
 problem.

\begin{problem} Classify or describe the Tridiagonal pairs.

\end{problem}

\noindent In \cite{TD00} we obtained the following  results on 
TD pairs.
 Referring to the  TD pair in Definition
\ref{def:tdrecall}, we showed $d=\delta$; we call this common
value the {\it diameter} of the pair. We showed that for
$0 \leq i \leq d$, the eigenspaces $V_i$ and $V^*_i$ have
the same dimension. Denoting this  common dimension by
$\rho_i$, 
 we showed 
the sequence $\rho_0, \rho_1, \ldots, \rho_d$ 
is symmetric and  unimodal; that is,
$\rho_i =\rho_{d-i}$ for  $0 \leq i \leq d$ and 
$\rho_{i-1} \leq \rho_{i}$ for $ 1 \leq i \leq d/2$.
We also obtained some results of a  representation theoretic nature, 
which we will discuss in the next section.

\section{Tridiagonal pairs in representation theory}

\medskip
\noindent In this section we survey how TD pairs arise
from irreducible finite dimensional modules
of the Lie algebra $sl_2$,  the Onsager algebra,
the quantum algebra $U_q(sl_2)$, and the $q$-Onsager algebra.
We define an algebra on two generators which we call the
{\it Tridiagonal algebra}. This algebra generalizes both
the Onsager and $q$-Onsager algebra. We show every
TD pair comes from an irreducible 
 finite dimensional 
module of a Tridiagonal algebra. 
For a more detailed discussion of the 
results in this section see \cite{TD00}.

\medskip
\noindent
We begin with $sl_2$.

\begin{example} \cite[Ex. 1.5]{TD00}
\label{ex:sl2easy} 
Assume $\K$ is algebraically closed with
characteristic 0, and let $L$ denote the Lie algebra
$sl_2(\K)$. Let $A$ and $A^*$ denote semi-simple 
elements in $L$ and assume $L$ is generated by
these elements. 
Let $V$ denote an irreducible finite dimensional
 $L$-module. Then
$A$ and $A^*$ act on $V$ as a Leonard  pair.

\end{example}

\noindent Referring to the previous example, it is not
hard to show  $A$ and $A^*$ must satisfy
\begin{eqnarray}
\lbrack A, \lbrack A, \lbrack A, A^* \rbrack \rbrack \rbrack &=&  b^2\lbrack
A, A^* \rbrack,
\label{eq:DG1}
\\
\lbrack A^*, \lbrack A^*, \lbrack A^*, A \rbrack \rbrack \rbrack &=& 
b^{*2}\lbrack
A^*, A \rbrack,
\label{eq:DG2}
\end{eqnarray}
where 
$\lbrack \;,\; \rbrack$ denotes  the Lie bracket,
and where $b$ and $b^*$ denote nonzero scalars in $\K$.
The relations
(\ref{eq:DG1}), 
(\ref{eq:DG2}) are known in statistical mechanics 
as the {\it Dolan-Grady relations}
 \cite{DateRoan2},
\cite{Dav},
\cite{Dolgra}, 
\cite{Roanmpi},
\cite{Ugl}.

\medskip
\noindent Using the Dolan-Grady relations, we obtain
TD pairs which are not necessarily Leonard pairs. 

\begin{example} \cite{TD00}
\label{def:onsager}
Assume $\K$ is algebraically closed 
with characteristic 0. Let $b$ and $b^*$ denote
nonzero scalars in $\K$.
Let $O$ denote the Lie algebra over $\K$ generated by
symbols $A$, $A^*$ subject to the Dolan-Grady relations
(\ref{eq:DG1}), 
(\ref{eq:DG2}). Let $V$ denote an irreducible finite dimensional
 $O$-module. Then $A, A^*$ act on $V$ as 
a TD pair.
\end{example}

\medskip
\noindent Referring to the above example,
the algebra $O$ is known in statistical mechanics as the 
Onsager algebra 
\cite{CKOns},
 \cite{DateRoan},
 \cite{DateRoan2},
\cite{Dav},
\cite{Per},
\cite{Roanmpi},
\cite{Ugl}. 
See 
 Davies
\cite{Dav}, Roan \cite{Roanmpi},
and Date and Roan   
\cite{DateRoan2}
for a detailed  description 
of the irreducible finite dimensional $O$-modules.

\medskip
\noindent Our next two examples are $q$-analogs
of the previous two. To set the stage,
we recall the  quantum algebra
$U_q(sl_2)$ and its modules.  For notational convenience, we 
will replace $q$ by $q^{1/2}$ and consider
$U_{q^{1/2}}(sl_2)$.

\begin{definition} \cite[p.122]{Kassel}
\label{def:uqsl2}
Assume  $\K$ is algebraically closed,
and let $q$ denote a 
nonzero scalar in $\K$ which is not a root of unity. 
Let $U_{q^{1/2}}(sl_2)$ denote the associative $\K$-algebra 
with 1 generated by symbols $ e,  f, k, k^{-1}$ subject
to the relations
\beast
kk^{-1} = k^{-1}k= 1,
\eeast
\beast
 ke = q ek,\qquad\qquad  kf = q^{-1}fk,
\eeast
\beast
ef - fe  = {{k-k^{-1}}\over {q^{1/2} - q^{-1/2}}}.
\eeast

\end{definition}
\noindent We now recall the irreducible finite dimensional modules
for 
 $U_{q^{1/2}}(sl_2)$.
\begin{lemma} \cite[p. 128]{Kassel}
\label{lem:uqmods}
With reference to Definition \ref{def:uqsl2},
there exists a family 
\begin{eqnarray}
V_{\varepsilon,d} \qquad \quad 
\varepsilon \in \lbrace 1,-1\rbrace, \qquad \quad  d = 0,1,2\ldots
\label{eq:uqmods}
\end{eqnarray}
of 
 irreducible finite dimensional 
 $U_{q^{1/2}}(sl_2)$
 modules with the following properties.
The module 
$V_{\varepsilon,d}$ has a basis $v_0, v_1, \ldots, v_d$ satisfying
$k v_i = \varepsilon q^{{d/2}-i}v_i $ for $0 \leq i \leq d$,
$f v_i = \lbrack i+1 \rbrack_q v_{i+1} $ for $0 \leq i \leq d-1$,
$f v_d = 0$,
$e v_i = \varepsilon \lbrack d-i+1 \rbrack_q v_{i-1} $ for $1 \leq i \leq d$,
$e v_0 = 0$, where 
\begin{equation}
\lbrack i \rbrack_q = 
{{q^{i/2} - q^{-i/2}}\over {q^{1/2}-q^{-1/2}}} \qquad 
\quad i \in \Z .
\label{eq:brackdef}
\end{equation}
Every irreducible finite dimensional module for 
 $U_{q^{1/2}}(sl_2)$ is isomorphic to exactly one of  
the modules (\ref{eq:uqmods}). (Referring to line
(\ref{eq:uqmods}), if 
$\K$ has characteristic 2
we interpret the set $\lbrace 1,-1 \rbrace $ as having a single element). 
\end{lemma}

\medskip
\noindent The  following result was proved by the present
author in \cite{Terint}
and is implicit in the results of  Koelink and Jan Der. Keugt
\cite{Koelink2}, \cite{Koelink4}.

\begin{example}
\cite{Koelink2}, \cite{Koelink4},
\cite{Terint} 
\label{ex:uqsl2} Referring to Definition
\ref{def:uqsl2} and Lemma 
\ref{lem:uqmods},
let  $\alpha, \alpha^*$ denote nonzero scalars in  $\K$, and 
set
\beast
A&=& \alpha f + {{ k}\over {q^{1/2}-q^{-1/2}}},\\
\eeast
\beast
A^*& =& \alpha^*e + {{k^{-1}}\over {q^{1/2}-q^{-1/2}}}.
\eeast
Let $d$ denote a nonnegative integer and pick
$\varepsilon \in \lbrace 1,-1 \rbrace $. 
 Then $A, A^*$ act on $V_{\varepsilon,d}$ as a Leonard pair
provided $\varepsilon \alpha \alpha^* $ is not among
$q^{(d-1)/2}, q^{(d-3)/2},\ldots, q^{(1-d)/2}$.
\end{example}
\noindent For related results concerning $U_{q^{1/2}}(sl_2)$ see 
\cite{Koelink3},
\cite{Koelink1},
\cite{koo3},
\cite[ch. 4]{Hjal}.

\medskip
\noindent Referring to Example  
\ref{ex:uqsl2}, one can show $A$ and $A^*$ satisfy
\begin{eqnarray} 
 A^3A^* - \lbrack 3\rbrack_q  A^2A^*A +  
 \lbrack 3\rbrack_q AA^*A^2  - A^*A^3   &=& 0,
\label{eq:qserre1} 
\\
 A^{*3}A -  \lbrack 3\rbrack_q  A^{*2}AA^* +  
 \lbrack 3\rbrack_q A^*AA^{*2}  - AA^{*3} &=& 0,
\label{eq:qserre2} 
\end{eqnarray}
where
$ \lbrack 3\rbrack_q= q+q^{-1}+1$. 
The equations 
(\ref{eq:qserre1}), 
(\ref{eq:qserre2}) 
are known as the {\it $q$-Serre  relations},
and are among the defining relations for the quantum affine algebra
$U_{q^{1/2}}({\widehat{sl}}_2)$. See  
\cite{CPqaa}, \cite{CPqaar} for more information about  this algebra.

\medskip
\noindent Using the $q$-Serre relations, we obtain TD pairs that are
not necessarily Leonard pairs.

\begin{example}\cite[Ex. 1.7]{TD00}
\label{ex:qOnsager}
Assume $\K$  is algebraically closed,
and let $q$ denote a nonzero element of $\K$ which is not a root of unity.
Let $O_q$ denote the associative $\K$-algebra with 1
generated by 
symbols $A$, $A^*$ subject to the $q$-Serre relations
(\ref{eq:qserre1}), 
(\ref{eq:qserre2}).
Let $V$ denote an irreducible finite dimensional module for  $O_q$,
and assume neither $A$ nor $A^*$ acts nilpotently on $V$.
Then $A, A^*$ act on $V$ as a TD pair.

\end{example}

\noindent Referring to the above example,
we call $O_q$  the {\it $q$-Onsager algebra}.

\medskip
\noindent 
The TD pairs that come from the above four examples
do not exhaust all TD pairs.
 To get all TD pairs,
we consider
 a pair of relations which
generalize both the 
Dolan-Grady relations and the
$q$-Serre relations. 
We call these the {\it Tridiagonal relations}.
In \cite{TD00} we proved the following result for an arbitrary field
$\K$.

\begin{theorem}\cite{TD00} 
\label{thm:gentd}
Let  $A, A^*$ denote a TD pair
over $\K$.
Then
there exists a sequence of  
 scalars
$\beta, \gamma, \gamma^*, \varrho, \varrho^* $ taken from $\K$
such that
\begin{eqnarray}
\lbrack A,A^2A^*-\beta AA^*A + A^*A^2 -\gamma (AA^*+A^*A)-\varrho A^*\rbrack 
&=&0,
\label{eq:TD1}
\\
\lbrack A^*,A^{*2}A-\beta A^*AA^* + AA^{*2} -
\gamma^* (A^*A+AA^*)-\varrho^* A \rbrack
&=&0,
\label{eq:TD2}
\end{eqnarray}
where $\lbrack r,s\rbrack$ means  $rs-sr$. 
The sequence 
$\beta, \gamma, \gamma^*, \varrho, \varrho^* $
is uniquely determined by the pair if the diameter is at least 3.
We refer to 
$\beta, \gamma, \gamma^*, \varrho, \varrho^* $ as a {\it parameter
sequence} for $A,A^*$. 
\end{theorem}

\noindent  
We call 
(\ref{eq:TD1}), 
(\ref{eq:TD2})  the 
{\it Tridiagonal relations} (or {\it TD relations}). 

\begin{remark}
\label{rem:dg}
 The Dolan-Grady relations (\ref{eq:DG1}), (\ref{eq:DG2}) are 
  the TD relations with parameters 
$\beta  = 2$, $\gamma = \gamma^*=0$,
$\varrho = b^2$, 
$\varrho^* = b^{*2}$, if we interpret the bracket 
 in (\ref{eq:DG1}), (\ref{eq:DG2}) as 
 $\lbrack r,s\rbrack = rs-sr$.   
The $q$-Serre relations 
(\ref{eq:qserre1}), 
(\ref{eq:qserre2}) 
are the TD relations with parameters
$\beta = q+q^{-1}$, $\gamma=\gamma^*=0$,
$\varrho=\varrho^*=0$.
\end{remark}

\begin{definition} 
\cite{TD00},
\cite{LS99}
\label{def:tdalg}
Let  
$\beta, \gamma, \gamma^*, \varrho, \varrho^* $ denote a sequence
of scalars taken from $\K$.
We let 
 $T=T(\beta, \gamma, \gamma^*,
\varrho, \varrho^*)$  denote the associative $\K$-algebra
with 1  generated by two symbols $A$, $A^*$
subject to the Tridiagonal relations 
(\ref{eq:TD1}), 
(\ref{eq:TD2}). We refer to $T$ as the {\it Tridiagonal algebra} 
(or {\it TD algebra})
over $\K$ with parameters 
$\beta, \gamma, \gamma^*, \varrho, \varrho^* $.
We refer to $A $ and $A^*$ as the {\it standard generators} of $T$.
\end{definition}

\noindent Reformulating Theorem
\ref{thm:gentd} using the above definition, we
can say the following.
Let $V$ denote a vector space over $\K$
with finite positive dimension. 
Let  $A, A^*$ denote a TD pair
on $V$ with parameter sequence
$\beta, \gamma, \gamma^*, \varrho, \varrho^*$.
By Theorem
\ref{thm:gentd}
the maps $A, A^*$ induce on  $V$ 
a module structure 
for the algebra 
$T(\beta, \gamma, \gamma^*, \varrho, \varrho^*)$.
This module is irreducible in view of 
Definition \ref{def:tdrecall}(iv).

\medskip
\noindent Given our above comments, the reader might wonder
if every irreducible finite dimensional  module for a 
TD algebra gives a TD pair. It turns out this is not
quite true. 
Our result is the following.

\begin{theorem}
\label{th:tmod}
Let  
 $\beta, \gamma, \gamma^*,
\varrho, \varrho^*$ denote scalars in $\K$, and 
assume $q$ is not a root of unity, where $q+q^{-1} = \beta$. 
Let $T$ denote the TD algebra over $\K$  
 with parameters $\beta, \gamma, \gamma^*,
\varrho, \varrho^*$ and standard generators $A, A^*$. 
Let $V$ denote an irreducible finite dimensional $T$-module
and assume each of $A, A^*$ is diagonalizable on $V$.
Then $A, A^*$ act on $V$ as a TD pair. 

\end{theorem}

\begin{proof} We show $A$ and $A^*$ satisfy conditions
(i)--(iv) of Definition \ref{def:tdrecall}. 
Condition (i) holds by assumption, so consider condition (ii).
Let $\Omega$ denote the set of distinct eigenvalues
of $A$ on $V$. For $\theta \in \Omega$, let $V_\theta$ denote the
corresponding eigenspace for $A$ and let $E_\theta$ denote the projection
of $V$ onto $V_\theta $. We recall $E_\theta - I$ vanishes on $V_\theta $,
where $I$ denotes the identity map on $V$, and that
$E_\theta $ vanishes on $V_\mu$ for all $\mu \in \Omega$,
$\mu \not= \theta $. 
Let $\theta, \mu$ denote distinct elements of $\Omega$. 
We  determine when  $E_\theta A^* E_\mu$ is zero.
To do this, we introduce the following
polynomial in two variables $x, y$.
\begin{equation}
p(x,y) = x^2-\beta x y +
y^2-\gamma (x+y)-\varrho.
\label{eq:pxy}
\end{equation}
We now assume $E_\theta A^* E_\mu \not=0$ and show $p(\theta, \mu) = 0$.
For notational convenience, set
\beast
C &=& 
A^2A^*-\beta AA^*A+A^*A^2 -\gamma(AA^*+A^*A) -\varrho A^*,
\eeast
and observe  $AC=CA$ by
(\ref{eq:TD1}). Using $E_\theta A = \theta E_\theta $  and 
 $AE_\mu  =  \mu E_\mu $, we have 
\beast
0 &=& E_\theta (AC-CA)E_\mu
\\
&=& (\theta-\mu)E_\theta CE_\mu,
\eeast
and since $\theta \not= \mu$,
\beast
0 &=& E_\theta CE_\mu
\\
&=&
E_\theta A^*E_\mu
(\theta^2-\beta\theta\mu +\mu^2 -
\gamma(\theta+\mu)-\varrho)
\\
&=&
E_\theta A^*E_\mu p(\theta,\mu).
\eeast
We assumed 
$E_\theta A^*E_\mu \not=0$,  so $p(\theta,\mu)=0$ as desired.
For all $\theta, \mu \in \Omega $, we define
 $\theta $ and  $\mu$ to be {\it adjacent} whenever 
$\theta \not=\mu$  and $p(\theta, \mu)=0$. Since
the polynomial $p$ is quadratic in its arguments,
each element in $\Omega $ is adjacent at most two
elements of $\Omega $. We show the adjacency relation
has no ``cycles''. Let $r$ denote an integer at least 3.
We claim there does not exist a sequence $\theta_1, \theta_2,\ldots,
\theta_r$ consisting  of distinct elements of $\Omega $, such
that $\theta_{i}$ is adjacent $\theta_{i+1}$ for $1 \leq i < r$
and $\theta_r$ is adjacent $\theta_1$. Suppose such a sequence
exists and consider the infinite sequence 
\begin{equation}
\theta_1, \theta_2,\ldots,
\theta_r,
\theta_1, \theta_2,\ldots,
\theta_r, \ldots
\label{eq:infseq}
\end{equation}
For $i = 0,1,2,\ldots $ let $\sigma_i$ denote the $i^{\hbox{th}}$
term in this sequence. For 
$i = 1,2,\ldots $ observe   $\sigma_{i-1}$  and $\sigma_{i+1} $
are distinct and adjacent $\sigma_i$, so they are the roots 
of $p(x,\sigma_i)$. Apparently $\sigma_{i-1} + \sigma_{i+1}$
is the opposite of the coefficient of $x$ in $p(x,\sigma_i)$.
Setting $y=\sigma_i$ in 
(\ref{eq:pxy}), we see 
this coefficient is $ -\beta \sigma_i -\gamma $, so
\begin{equation}
\sigma_{i-1} - \beta \sigma_i + \sigma_{i+1} = \gamma
\qquad \qquad i= 1,2,\ldots
\label{eq:sigrec}
\end{equation}
By definition $\beta = q + q^{-1}$, and we assume 
$q \not\in \lbrace 1,-1\rbrace $,
so $\beta \not\in \lbrace 2,-2\rbrace $.
Solving the linear recurrence 
(\ref{eq:sigrec}) 
using this, we find there exists
scalars $a,b,c$ in the algebraic closure of $\K$ such that
\begin{equation}
\sigma_i = a + bq^i + cq^{-i} \qquad \qquad i = 0,1,2,\ldots
\label{eq:sigcl}
\end{equation}
We show 
(\ref{eq:sigcl}) is inconsistent with the periodic nature of
(\ref{eq:infseq}). For $i=0,1,2,\ldots $ we have 
$\sigma_i = \sigma_{i+r}$. Using 
(\ref{eq:sigcl}) we find
$\sigma_i - \sigma_{i+r}$ equals $1-q^r$ times
\begin{equation}
bq^i - cq^{-i-r}.
\label{eq:nocycle}
\end{equation}
We assume $q^r\not=1$, so  the expression
(\ref{eq:nocycle}) is zero. Setting  the expression
(\ref{eq:nocycle}) equal to  zero 
for $i=0,1$ and using $q\not=1$, $q\not=-1$, we routinely find
$b=0, c=0$. But this is inconsistent with 
(\ref{eq:sigcl}) and the requirement $\sigma_0 \not=\sigma_1$.
We have now proved the claim.
From our above remarks, there exists an ordering 
$\theta_0, \theta_1, \ldots, \theta_d $ of all the elements of
$\Omega $ such that $\theta_i, \theta_j$ are not adjacent for
$|i-j|>1 $, $(0 \leq i,j\leq d)$. From our preliminary comments
and abbreviating $E_i $ for $E_{\theta_i}$,
\beast
E_iA^*E_j = 0 \quad \hbox{if} \quad  |i-j|>1, \qquad \quad (0 \leq i,j\leq d).
\eeast
From the above line and abbreviating $V_i$ for $V_{\theta_i}$, we obtain
\beast
A^* V_i \subseteq V_{i-1} + V_i + V_{i+1} \qquad \qquad (0 \leq i \leq d),
\eeast
where $V_{-1} = 0$, $V_{d+1}=0$. We now see $A, A^*$ satisfy
Definition \ref{def:tdrecall}(ii). 
Interchanging the roles of $A$ and $A^*$ in the above argument, we see 
 $A, A^*$ satisfy
Definition \ref{def:tdrecall}(iii).     
Observe $A, A^*$ satisfy
Definition \ref{def:tdrecall}(iv), since we assume
$V$ is irreducible as a $T$-module.
We have now shown $A, A^*$ satisfy conditions (i)--(iv) of
Definition  \ref{def:tdrecall}, so $A,A^*$ act on $V$ as a TD pair.

\end{proof}

\noindent We finish this section with a few comments on the 
Tridiagonal relations 
(\ref{eq:TD1}),
(\ref{eq:TD2}). 
These relations previously appeared in \cite{TersubIII}.
In that paper, the author
considers a combinatorial  object called  
a $P$- and $Q$-polynomial association scheme
\cite{BanIto}, \cite{bcn},
\cite{Leopandq},
\cite{Tercharpq},
\cite{Ternew}.
He shows 
 that for these schemes the adjacency matrix $A$ 
 and a certain diagonal matrix  $A^*$
satisfy
(\ref{eq:TD1}), 
(\ref{eq:TD2}). 
In this context the algebra 
generated by $A$ and $A^*$ is known as the subconstituent
algebra or the Terwilliger algebra
\cite{TersubI},
\cite{TersubII},
\cite{TersubIII}. See 
\cite{Cau},
\cite{Col},
\cite{CurNom},
\cite{Curbip1},
\cite{Curbip2},
\cite{Curspin},
\cite{Go},
\cite{HobIto},
\cite{Tan} for more information on this algbra.

\medskip
\noindent 
We  mention the relations
(\ref{eq:TD1}),
(\ref{eq:TD2})  are satisfied by the
generators of 
both the classical
and quantum Quadratic  Askey-Wilson algebra
introduced by   Granovskii, Lutzenko, and 
Zhedanov \cite{GYLZmut}.
See 
\cite{GYZnature},
\cite{GYZTwisted},
\cite{GYZlinear},
\cite{GYZspherical},
\cite{Zhidd}, 
\cite{ZheCart},
\cite{Zhidden}
for more information on this algebra.

\medskip
\noindent 
The relation (\ref{eq:TD2}) previously appeared in
the work of Grunbaum and Haine on the ``bispectral problem''
\cite{GH7},
\cite{GH6}.
See
\cite{GH4},
\cite{GH5},
\cite{GH1}, 
\cite{GH3},
\cite{GH2} 
for related work.

\section{The parameter sequence and the eigenvalues}

\medskip
\noindent In this section we obtain the following
results concerning the  
parameter sequence of a  TD pair.
We consider 
a transformation on TD pairs and determine the effect on
the parameter sequence. We recall the eigenvalue and
dual eigenvalue sequences of a TD pair, and show
these sequences satisfy a three term recurrence.
We show how to get the parameter sequence of a TD pair from its
eigenvalue and dual eigenvalue sequences. 
We characterize the TD pairs satisfying the Dolan-Grady relations 
in terms of their eigenvalue sequence and dual eigenvalue sequence.
We obtain a similar characterization of the TD pairs 
satisfying the $q$-Serre relations.

\medskip
\noindent
 We begin with an observation. Let $V$ denote 
a vector space over $\K$ with finite positive dimension,
and let $A,A^*$ denote a TD pair on $V$. Let 
$r,s, r^*,  s^* $ denote scalars in $\K$ with each of $r, r^*$ nonzero.
Then the ordered pair
\begin{equation}
rA+sI, \qquad  r^*A^*+s^*I
\label{eq:reldef}
\end{equation}
is a TD pair on $V$.
The above transformation 
 has the following effect on the parameter 
sequence. 
Let $\beta, \gamma, \gamma^*,
\varrho, \varrho^*$ 
denote a parameter sequence for $A,A^*$,
as in 
Theorem \ref{thm:gentd}.
Then the following is a parameter sequence for 
the TD pair   
(\ref{eq:reldef}).
\beast
&&\beta, \qquad 
 r \gamma + s(2-\beta), 
\qquad 
 r^* \gamma^* + s^*(2-\beta), 
\\
&& r^2\varrho-2rs\gamma+s^2(\beta-2), \qquad 
r^{*2}\varrho^*-2r^*s^*\gamma^*+s^{*2}(\beta-2). 
\eeast
A given parameter sequence can often be simplified by
transforming it in the above fashion.
How far one can go in this direction depends on the field $\K$,
so we will not give all the details. Instead, we illustrate with an example.
 Let us say a parameter sequence 
$\beta, \gamma, \gamma^*, \varrho, \varrho^*$ is {\it reduced}
whenever $\gamma = \gamma^*=0$.
Let $A,A^*$ denote a TD pair on $V$
 with parameter sequence 
 $\beta, \gamma, \gamma^*, \varrho, \varrho^*$.
Then the  TD pair below
has a reduced parameter sequence provided $\beta \not=2$.
\beast
A + \gamma (\beta - 2)^{-1}I, \qquad \qquad 
A^* + \gamma^* (\beta - 2)^{-1}I.
\eeast

\noindent We now consider the relationship between 
the parameter sequence
and the eigenvalues of $A$ and $A^*$. 

\begin{definition} \cite{TD00}
\label{def:eigvals}  
Let $A,A^*$ denote a TD pair with diameter $d$.
By an {\it eigenvalue sequence} for $A,A^*$, we mean an ordering
 $\theta_0, \theta_1, \ldots, \theta_d$ of    
the distinct eigenvalues of $A$   
that satisfies Definition
\ref{def:tdrecall}(ii), where for $0 \leq i \leq d$ 
the $V_i$ in that definition
denotes the eigenspace of
$A$ associated with $\theta_i$.
We remark that if 
 $\theta_0, \theta_1, \ldots, \theta_d$ is an eigenvalue
 sequence for $A,A^*$ then so is 
 $\theta_d, \theta_{d-1}, \ldots, \theta_0$,
 and $A,A^*$ has no further eigenvalue
 sequence. By a {\it dual eigenvalue sequence} for $A,A^*$, we mean
 an eigenvalue sequence for  $A^*, A$.

\end{definition}

\medskip  
\noindent 
In \cite{TD00} we obtained
the following  result.

\begin{theorem} \cite{TD00}
\label{thm:eigprel}
let $A,A^*$ denote a TD pair,  
with eigenvalue sequence $\theta_0, \theta_1, \ldots, \theta_d$ and
dual eigenvalue sequence
 $\theta^*_0, \theta^*_1, \ldots, \theta^*_d$.
Then the expressions
\begin{equation}
{{\theta_{i-2}-\theta_{i+1}}\over {\theta_{i-1}-\theta_i}},\qquad \qquad  
 {{\theta^*_{i-2}-\theta^*_{i+1}}\over {\theta^*_{i-1}-\theta^*_i}} 
 \qquad  \qquad 
\label{eq:expsbetaplusone}
\end{equation} 
 are equal and independent of $i$ for $\;2\leq i \leq d-1$.  

\end{theorem}

\noindent 
In the following theorem, we explain how to obtain a parameter
sequence for a given TD pair from its eigenvalue sequence and
dual eigenvalue sequence. 
To prepare for this result we make some comments.
For the moment, let
$\theta_0, \theta_1, \ldots, \theta_d$ denote any finite sequence of
distinct scalars in $\K $, and  assume    
\begin{equation}
{{\theta_{i-2}-\theta_{i+1}}\over {\theta_{i-1}-\theta_i}}
\label{eq:arbseq}
\end{equation}
is independent of $i$ for $2 \leq i \leq d-1$. Denoting by $\beta +1$
the common value of the expressions
(\ref{eq:arbseq}), we find 
$\theta_{i-1}-\beta \theta_i + \theta_{i+1}$ is independent of
$i$ for $1 \leq i \leq d-1$. Denoting this common value
by $\gamma $, we have
\beast
\gamma = \theta_{i-1}-\beta \theta_i + \theta_{i+1} \qquad \qquad (1 \leq i \leq d-1).
\eeast
We claim there exists a scalar $\varrho \in  \K $ such that
\begin{equation}
\varrho=
\theta^2_{i-1}-\beta \theta_{i-1}\theta_i+\theta_i^2-
\gamma (\theta_{i-1}+\theta_i) \qquad \qquad (1 \leq i \leq d).
\label{eq:vrhoarb}
\end{equation}
To  see this, let $p_i$ denote the expression on the right in   
(\ref{eq:vrhoarb}), and observe
\begin{equation}
p_i - p_{i+1} = (\theta_{i-1}-\theta_{i+1})(\theta_{i-1}-\beta \theta_i+\theta_{i+1} - \gamma)
\label{eq:pdiff}
\end{equation}
for $1 \leq i\leq d-1$. From 
(\ref{eq:pdiff}) we see $p_i$ is independent of $i$ for $1 \leq i \leq d$.
Denoting the common value of the $p_i$ by $\varrho$, we obtain
(\ref{eq:vrhoarb}).
With these comments in mind, we present the next result.

\begin{theorem} \cite{TD00}
\label{thm:eiginfo}
Let $A,A^*$ denote a TD pair over  
 $\K$,
with eigenvalue sequence $\theta_0, \theta_1, \ldots, \theta_d$
and dual eigenvalue sequence 
 $\theta^*_0, \theta^*_1, \ldots, \theta^*_d$.
Let 
 $\beta, \gamma, \gamma^*,
\varrho, \varrho^*$ 
denote a sequence of  scalars taken from $\K$.
Then this is a parameter sequence for $A,A^*$ 
 if and only if
(i)--(v) hold below.
\begin{enumerate}
\item
The expressions
\beast
{{\theta_{i-2}-\theta_{i+1}}\over {\theta_{i-1}-\theta_i}},\qquad \qquad  
 {{\theta^*_{i-2}-\theta^*_{i+1}}\over {\theta^*_{i-1}-\theta^*_i}} 
 \qquad  \qquad 
\eeast
both equal $\beta +1$ for $\;2\leq i \leq d-1$.  
\item
$\gamma = \theta_{i-1}-\beta \theta_i + \theta_{i+1} \qquad \qquad (1 \leq i \leq d-1)$,
\item
$\gamma^* =
\theta^*_{i-1}-\beta \theta^*_i + \theta^*_{i+1} 
\qquad \qquad (1 \leq i \leq d-1),$
\item
$\varrho= \theta^2_{i-1}-\beta \theta_{i-1}\theta_i+\theta_i^2-\gamma (\theta_{i-1}+\theta_i) \qquad \qquad (1 \leq i \leq d),$
\item
$\varrho^*= \theta^{*2}_{i-1}-\beta \theta^*_{i-1}\theta^*_i+\theta_i^{*2}-
\gamma^* (\theta^*_{i-1}+\theta^*_i) \qquad \qquad (1 \leq i \leq d). $
\end{enumerate}

\end{theorem}

\noindent To get another point of view, we consider
the eigenvalues and dual eigenvalues in parametric form.
Solving the linear recurrence in 
Theorem \ref{thm:eigprel}, we obtain the following.


\begin{theorem}
\cite{TD00}
\label{thm:eigclosed}
Let $A,A^*$ denote a TD pair over $\K$, with 
 eigenvalue sequence  
$\theta_0, \theta_1,\ldots, \theta_d$ and dual
eigenvalue sequence 
$\theta^*_0, \theta^*_1,\ldots, \theta^*_d$.
Then these sequences are given by 
Case I, II, or III below.

\medskip
\noindent Case I: 
\beast
\theta_i &=& a + bq^i + cq^{-i},
\qquad \qquad \qquad q\not=0, \quad  q \not= 1, \quad  q\not=-1\\
\theta^*_i &=& a^* + b^*q^i + c^*q^{-i}. 
\eeast
\noindent Case II:  
\beast
\theta_i &=& a + bi + ci(i-1)/2,  
\\
\theta^*_i &=& a^* + b^*i + c^*i(i-1)/2.  
\eeast
\noindent Case III:  The characteristic of $\K$ is not 2,  and 
\beast
\theta_i &=& a + b(-1)^i + ci(-1)^i,  
\\
\theta^*_i &=& a^* + b^*(-1)^i + c^*i(-1)^i.  
\eeast

\noindent In the above formulae, the scalars $q$ and 
$a,b,c, a^*, b^*, c^*$ are in the algebraic closure 
of $\K$.  Concerning Case II, if the characteristic of $\K$ equals 2,  
 we interpret the expression $i(i-1)/2 $ as 0 if $i=0$ or $i=1$ (mod 4),
and as 1 
 if $i=2$ or $i=3$ (mod 4).
\end{theorem}

\noindent
Evaluating the data in 
Theorem \ref{thm:eiginfo} using
Theorem \ref{thm:eigclosed}, we obtain the following.

\begin{lemma}
\label{cor:abc}  
Referring to Theorem
\ref{thm:eigclosed}, the sequence  
 $\beta, \gamma, \gamma^*,
\varrho, \varrho^*$ given below is  a 
parameter sequence 
for $A,A^*$.

\medskip
\noindent Case I:
\beast
\beta &=& q+q^{-1},
\\
\gamma &=& -a (q-1)^2q^{-1}, 
\\
\varrho &=& a^2(q-1)^2q^{-1}- bc(q-q^{-1})^2. 
\eeast
\noindent Case II:
\beast
\beta &=& 2, 
\\
\gamma &=&  c,
\\
\varrho &=& b^2-bc-2ac.
\eeast
\noindent Case III:
\beast
\beta &=& -2, 
\\
\gamma &=& 4a,
\\
\varrho &=& c^2-4a^2.
\eeast
To get $\gamma^*$, $\varrho^*$, replace $a,b,c$ in the above lines
by $a^*,b^*,c^*$.
\end{lemma}

\noindent At the beginning of this section we defined the notion
of a reduced parameter sequence. Concerning this we have the following
result.

\begin{corollary}
\label{lem:redmeaning}
The parameter sequence given in 
Lemma \ref{cor:abc} is reduced if and only if the following
hold.

\medskip
\noindent Cases I, III: \qquad $a = 0$, $\;a^*=0$.

\medskip
\noindent Case II: \qquad $c = 0$, $\;c^*=0$.
%
\end{corollary}

\noindent 
We finish this section with some comments on the
 Dolan-Grady relations
and the
$q$-Serre relations.
For the moment let
$\theta_0, \theta_1, \ldots, \theta_d$ denote any finite sequence
of scalars in $\K$.
Let $b$ denote a scalar in $\K$. We say 
the sequence $\theta_0, \theta_1, \ldots, \theta_d$ is in {\it $b$-arithmetic
progression}  whenever  $\theta_i = \theta_{i-1}  + b$ 
for $1 \leq i \leq d$.
%
Let $q$ denote a nonzero scalar in $\K$.
We say  
$\theta_0, \theta_1, \ldots, \theta_d$ is in {\it $q$-geometric progression} 
whenever $\theta_i =  \theta_{i-1} q$ for 
$1 \leq i   \leq d$.
%
%
%
\begin{lemma}
\label{lem:dgcom}
Let $A,A^*$ denote a TD pair over $\K$ and let 
 $b, b^*$ denote nonzero scalars in $\K$. Then (i), (ii)
 below are equivalent.
\begin{enumerate}
\item 
 $A,A^*$ satisfy the Dolan-Grady relations (\ref{eq:DG1}), 
(\ref{eq:DG2}) for $b$ and $b^*$. 
\item There exists an eigenvalue sequence for $A,A^*$ 
which is in 
$b$-arithmetic progression, and there exists a dual eigenvalue sequence 
 for $A,A^*$  which is 
in
$b^*$-arithmetic progression.
\end{enumerate}
\end{lemma}

\begin{proof} $(i)\rightarrow (ii) $ By Remark
\ref{rem:dg} we see
$2, 0,0,b^2, b^{*2}$ is a parameter sequence for $A,A^*$.
Applying Theorem
\ref{thm:eiginfo} to this sequence, we routinely obtain 
the result.

\noindent $(ii)\rightarrow (i)$. Setting $c=0$, $c^*=0$ in
Case II of Lemma
\ref{cor:abc}, we find $2,0,0,b^2,b^{*2}$ is a parameter
sequence for $A,A^*$. The result now follows in view
of Remark 
\ref{rem:dg}.

\end{proof}

\begin{lemma}
\label{lem:qscom}
Let $A,A^*$ denote a TD pair over $\K$ and let
 $q$ denote a nonzero scalar in $\K$ other than
 $1,-1$.
Then (i), (ii) below are equivalent. 
\begin{enumerate}
\item
$A,A^*$ satisfy the
$q$-Serre relations
(\ref{eq:qserre1}), 
(\ref{eq:qserre2}).
\item There exists an 
 eigenvalue sequence for
$A,A^*$ which is  
 in $q$-geometric progression,
 and
there exists a dual 
 eigenvalue sequence for
$A,A^*$ which is  
 in $q$-geometric progression.
\end{enumerate}

\end{lemma}

\begin{proof} 
Similar to the proof of Lemma  \ref{lem:dgcom}.

\end{proof}

%
%
%
%
%
%
%

\section{Infinite dimensional modules for Tridiagonal algebras}

\medskip
\noindent 
In this section we consider  infinite dimensional
irreducible modules for 
TD algebras. These modules seem more complicated than the finite
dimensional ones, so we do not attempt a general theory. Instead,
we motivate further study by presenting two examples obtained
from certain orthogonal polynomials of the Askey scheme. Our first
example, based on the Hermite polynomials, is very simple and
meant to illustrate the ideas involved. Our second and main
example is based on the Askey-Wilson polynomials.
We would like to acknowledge that the results of this section
can be readily obtained from  any of the following
works:
 Granovskii, Lutzenko, and 
Zhedanov \cite{GYLZmut},
Grunbaum and Haine
\cite{GH7},
Noumi and Stokman \cite{NStok}. Our goal here is to emphasize the
role played by the Tridiagonal relations.

\medskip
\noindent Throughout this section $x$ will denote an indeterminant
and $\K\lbrack x \rbrack $ will denote the $\K$-algebra consisting
of all polynomials in $x$ that have coefficients in $\K$.

\medskip
\noindent {\bf Example 5.1}. {\bf The Hermite polynomials}.
For this example assume $\K$
has characteristic 0. Let $H_0, H_1, \ldots $ denote the polynomials
in $\K \lbrack x \rbrack $ satisfying 
\begin{equation}
x H_n = H_{n+1} + 2n H_{n-1}, \qquad   \quad n = 0,1,2,\ldots
\label{eq:hrec}
\end{equation}
where $H_0=1$, $H_{-1} = 0$. We refer to $H_n$ as the $n^{\hbox{th}}$
{\it Hermite polynomial} \cite{KoeSwa}. The first four Hermite polynomials are 
\beast
H_0=1,\qquad H_1 = x, \qquad H_2 = x^2-2, \qquad H_3=x^3-6x.
\eeast
The Hermite polynomials satisfy a second order differential equation,
obtained as follows. Let $D: \K\lbrack x \rbrack \rightarrow 
\K\lbrack x \rbrack  $ denote the derivative map. We recall
$D$ is a linear transformation satisfying $D x^n = nx^{n-1} $ for
$n = 0,1,2,\ldots $ It is well known
\begin{equation}
D H_n = n H_{n-1} \qquad \qquad n=0,1,2\ldots 
\label{eq:dh}
\end{equation}
This can be proved by induction on $n$, or see \cite[p. 280]{AAR}. Combining
(\ref{eq:hrec}), 
(\ref{eq:dh}), we find
\begin{equation}
(x-2D) H_n =  H_{n+1} \qquad \qquad n=0,1,2\ldots 
\label{eq:uph}
\end{equation}
Combining 
(\ref{eq:dh}), 
(\ref{eq:uph}), we obtain
\begin{equation}
(x-2D)D H_n =  nH_n \qquad \qquad n=0,1,2\ldots 
\label{eq:hdiff}
\end{equation}
This is the differential equation satisfied by the Hermite polynomials.
We now endow the vector space $\K\lbrack x \rbrack $ with a module
structure for a certain TD algebra. Let $A$ and $A^*$ denote 
the following linear transformations.
\beast
\K\lbrack x \rbrack \qquad &\rightarrow &   \qquad  \quad
\K\lbrack x \rbrack   
\\
A: \quad \qquad f \;\; \qquad  &\rightarrow &\qquad \quad\; xf,
\\
A^*: \quad \qquad f\;\; \qquad & \rightarrow  &\quad (x-2D)Df.
\eeast
With respect to the basis $H_0, H_1, \ldots$ the matrices
representing $A$ and $A^*$ are, respectively,
\begin{equation}
\left(
\begin{array}{ c c c c c c}
0 & 2  &      &      &   & \\
1  & 0  &  4   &      &   &  \\
  & 1  &  0    & 6  &   & \\
  &   &  1     & \cdot  & \cdot   & \\
  &   &           &  \cdot & \cdot &  \\
 &   &   &   &  &   
\end{array}
\right),
\qquad \hbox{diag}(0, 1, 2, \ldots ).
\label{eq:hmat}
\end{equation}
These are obtained from
(\ref{eq:hrec})  and
(\ref{eq:hdiff}). Representing $A$ and $A^*$ by the 
matrices in 
(\ref{eq:hmat}),
 we routinely find
\begin{eqnarray}
A^2A^*- 2 AA^*A + A^*A^2 &=& -4I,
\label{eq:htd1}
\\
A^{*2}A- 2 A^*AA^* + AA^{*2} -A&=& 0,
\label{eq:htd2}
\end{eqnarray}
where $I$ denotes the identity map on 
$\K\lbrack x \rbrack   $. Observe the left side of
(\ref{eq:htd1}) commutes with $A$ and the left side of
(\ref{eq:htd2}) commutes with $A^*$. Therefore $A$ and $A^*$
satisfy the TD relations 
(\ref{eq:TD1}), (\ref{eq:TD2}) for $\beta = 2$,
$\gamma=\gamma^*=0$, $\delta=0$, $\delta^*=1$. Apparently
$A$ and $A^*$ induce a module structure 
on $\K\lbrack x \rbrack $ for the  algebra $T(2,0,0,0,1)$. This
module is irreducible in view of the following lemma.
The proof of this lemma  is routine and omitted. 

\medskip
\noindent {\bf Lemma 5.2}
{\it Let $V$ denote an infinite dimensional vector space over $\K$.  
Let $A:V \rightarrow V$ and $A^*:V \rightarrow V$ denote linear
transformations. Assume $V$ has a basis $v_0, v_1,\ldots $ with
respect to which (i) the matrix representing $A$ is irreducible
tridiagonal, and (ii) the matrix representing $A^*$ is diagonal,
with diagonal entries mutually distinct. Then there does not
exist a subspace $W \subseteq V$ such that $AW\subseteq W$ and
$A^*W\subseteq W$, other than $W=0$ and $W=V$.
}

\medskip
\noindent {\bf Example 5.3}. {\bf The Askey-Wilson polynomials}. 
For this example assume 
the field $\K$ is arbitrary. Let $q, a,b,c,d$ denote nonzero scalars
in $\K$. To avoid degenerate situations, we assume $q$ is not
a root of unity, and that none of $ab,ac,ad,bc,bd,cd,abcd$ is an integral
power of $q$. For $n=0,1,2\ldots $ let $p_n = p_n(x;a,b,c,d)$ denote
the polynomial in $\K\lbrack x \rbrack $ given by
\beast
p_n=
 {{}_4}\phi_3\Biggl({{q^{-n}, abcdq^{n-1},ay,ay^{-1}}\atop {ab,ac,ad}} \bigg|\; q;q\Biggr),
\eeast
where $x=y+y^{-1}$. We follow the standard notation for basic hypergeometric
series given in \cite{gasperrahmanbk}. We refer to $p_n$ as the $n^{\hbox{th}}$ {\it Askey-Wilson
polynomial} with parameters $a,b,c,d$ 
\cite{AskAW},
\cite{gasperrahmanbk},
\cite{KoeSwa}.
The first three Askey-Wilson polynomials are $p_0=1$,
\beast
&&p_1 \; =\; 1 - {{(1-abcd)(1-ax+a^2)}\over {(1-ab)(1-ac)(1-ad)}},
\\
&&p_2 \;=\; 1 - 
 {{(1+q^{-1})(1-abcdq)(1-ax+a^2)}\over {(1-ab)(1-ac)(1-ad)}}
\\
&& \qquad \qquad \qquad + \; 
 {{(1-abcdq)(1-abcdq^2)(1-ax+a^2)(1-axq+a^2q^2)}\over 
 {q(1-ab)(1-abq)(1-ac)(1-acq)(1-ad)(1-adq)}}.
\eeast
\noindent The Askey-Wilson polynomials satisfy the following three
term recurrence \cite{AskAW}. For $n=0,1,2,\ldots $ we have
\begin{equation}
x p_n = b_n p_{n+1} + a_n p_n + c_n p_{n-1},
\label{eq:awrec}
\end{equation}
where $p_{-1}=0$, and where
\begin{eqnarray}
b_n &=& {{(1-abq^n)(1-acq^n)(1-adq^n)(1-abcdq^{n-1})}
\over {a(1-abcdq^{2n-1})(1-abcdq^{2n})}},
\label{eq:awbn}
\\
c_n &=& {{a(1-q^n)(1-bcq^{n-1})(1-bdq^{n-1})(1-cdq^{n-1})}
\over {(1-abcdq^{2n-2})(1-abcdq^{2n-1})}},
\label{eq:awcn}
\\
a_n &=& a+a^{-1} -b_n -c_n.
\label{eq:awan}
\end{eqnarray}
We observe the denominators in 
(\ref{eq:awbn}), 
(\ref{eq:awcn}) are not zero, 
and that $b_{n-1}c_n\not=0$ for $n=1,2,\ldots $
The Askey-Wilson polynomials satisfy a $q$-difference equation, obtained
as follows. Let $y$ denote an indeterminant, and let 
$\K\lbrack y,y^{-1} \rbrack $ denote the $\K$-algebra consisting of all
Laurent polynomials in $y$ that have coefficients in $\K$. Identifying
$x$ with $y+y^{-1}$, we view 
$\K\lbrack x\rbrack $ 
as the subalgebra of 
$\K\lbrack y,y^{-1} \rbrack $
generated by $x$. From this point of view 
$\K\lbrack x\rbrack $  has basis 
\beast
1, \quad y+y^{-1}, \quad y^2+y^{-2}, \quad y^3+y^{-3}, \quad \ldots
\eeast
Let $\tau$ denote the $\K$-algebra automorphism of 
$\K\lbrack y,y^{-1} \rbrack $ satisfying $\tau (y) = qy$.
Let 
$\D$ denote the restriction of the following map to 
$\K\lbrack x\rbrack $:
\begin{equation}
\phi(y)(\tau - I) \;+ \;\phi(y^{-1})(\tau^{-1}-I)\;+\;(1+abcdq^{-1})I.
\label{eq:dddef}
\end{equation}
Here $I$ denotes the identity map and 
\beast
\phi(y) = {{(1-ay)(1-by)(1-cy)(1-dy)}\over {(1-y^2)(1-qy^2)}}.
\eeast
In concrete terms, $\D :  
\K\lbrack x\rbrack \rightarrow  
\K\lbrack x\rbrack $ is the linear transformation satisfying
$\D (1) = 1+abcdq^{-1}$, and for $n=1,2,3 \ldots  $
\beast
&&\D(y^n+y^{-n}) = (1+abcdq^{-1})(y^n+y^{-n})
\\
&& \qquad \qquad 
+\;{{(1-ay)(1-by)(1-cy)(1-dy)}\over {(1-y^2)(1-qy^2)}}(q^ny^n-y^n+q^{-n}y^{-n}-y^{-n})
\\
&& \qquad \qquad 
+\;{{(1-ay^{-1})(1-by^{-1})(1-cy^{-1})(1-dy^{-1})}\over
{(1-y^{-2})(1-qy^{-2})}}(q^{-n}y^n-y^n+q^{n}y^{-n}-y^{-n}).
\eeast
 The map $\D$ is known as the {\it Askey-Wilson
second order $q$-difference operator} with parameters $a,b,c,d$
\cite{AskAW}.
The Askey-Wilson polynomials  mentioned above
are the eigenvectors for $\D$. In fact,
for $n=0,1,2,\ldots $ we have
\begin{equation}
\D \,p_n = \theta^*_n\, p_n,
\label{eq:awpolyev}
\end{equation}
where 
\begin{equation}
\theta^*_n = q^{-n} + abcdq^{n-1}.
\label{eq:thsaw}
\end{equation}
For a proof of this see \cite{AskAW}. Line 
(\ref{eq:awpolyev}) is the above-mentioned $q$-difference equation
satisfied by the Askey-Wilson polynomials. We now endow the vector
space 
$\K\lbrack x\rbrack $ with the module structure of a certain
TD algebra. Let $A$ and $A^*$ denote the following linear
transformations.
\beast
\K\lbrack x \rbrack \qquad &\rightarrow &   \qquad 
\K\lbrack x \rbrack   
\\
A: \quad \qquad f\;\; \qquad  &\rightarrow & \qquad \;xf,
\\
A^*: \quad \qquad f\;\; \qquad & \rightarrow  &\qquad \;\D f.
\eeast
We show $A$ and $A^*$ satisfy the TD relations
\begin{eqnarray}
\lbrack A,A^2A^*-(q+q^{-1}) AA^*A + A^*A^2 +(q-q^{-1})^2 A^*\rbrack 
&=&0,
\label{eq:awTD1}
\\
\lbrack A^*,A^{*2}A-(q+q^{-1}) A^*AA^* + AA^{*2} 
+abcdq^{-1}(q-q^{-1})^2 A \rbrack
&=&0,
\label{eq:awTD2}
\end{eqnarray}
where we recall
$\lbrack r,s\rbrack $ means $rs-sr$.
We begin with 
(\ref{eq:awTD2}). Observe $p_0, p_1, \ldots $ is a basis for the
vector space 
$\K\lbrack x \rbrack $. With respect to this basis, the matrices
representing $A$ and $A^*$ are 
\begin{equation}
\left(
\begin{array}{ c c c c c c}
a_0 & c_1  &      &      &   & \\
b_0 & a_1  &  c_2   &      &   &  \\
  & b_1  &  a_2    & \cdot  &   & \\
  &   &  \cdot    & \cdot  & \cdot   & \\
  &   &           &  \cdot & \cdot &  \\
 &   &   &   &  &   
\end{array}
\right),
\qquad \qquad \hbox{diag}(\theta^*_0, \theta^*_1, \theta^*_2, \ldots ),
\label{eq:awmat}
\end{equation}
respectively,
where the $a_n, b_n, c_n$ are from
(\ref{eq:awbn})--(\ref{eq:awan}), and where the $\theta^*_n$ are from
(\ref{eq:thsaw}). Let $\Delta $ denote the matrix representing
the left side of 
(\ref{eq:awTD2}) with respect to the basis $p_0, p_1,\ldots $ We show
$\Delta = 0$. To do this, we show each entry of
$\Delta $ is zero.
Using the matrices 
(\ref{eq:awmat})
 we find that  for nonnegative integers $r,s$,  
 the matrix $\Delta $ has $r,s$ entry
\begin{equation}
\theta^{*2}_r-(q+q^{-1})\theta^*_r\theta^*_s
+ \theta^{*2}_s+abcdq^{-1}(q-q^{-1})^2
\label{eq:dent}
\end{equation}
times
$(\theta^*_r-\theta^*_s)
B_{rs}$,
where $B$ denotes the  matrix on the left in
(\ref{eq:awmat}). 
Observe $B$ is tridiagonal, so $B_{rs} = 0$ for $|r-s|>1$.
In view of 
(\ref{eq:thsaw}) the expression 
(\ref{eq:dent}) is zero 
for $|r-s|=1$. Of course
$\theta^*_r-\theta^*_s$  is zero  for $r=s$.
Apparently $\Delta $ has all entries  zero, so $\Delta = 0$, and
(\ref{eq:awTD2}) follows.
We now show 
(\ref{eq:awTD1}).
It is possible to verify  
(\ref{eq:awTD1}) using the above matrix representations
of $A$ and $A^*$, but the calculation is tedious.
Instead, we  follow the argument of   
 Grunbaum and Haine
\cite{GH7}
 who obtained 
(\ref{eq:awTD1}) in a different context.
 To obtain 
(\ref{eq:awTD1}), we show $A$ commutes with
\begin{equation}
A^2A^*-(q+q^{-1}) AA^*A + A^*A^2 +(q-q^{-1})^2 A^*.
\label{eq:acom}
\end{equation}
To do this, we recall the map
$A$ represents multiplication by $x$, and
$A^*$ represents the Askey-Wilson operator
$\D$ given in 
(\ref{eq:dddef}). Consider the terms in 
(\ref{eq:dddef}).  
Recall the map $\tau$ in that line
is the automorphism of 
$\K\lbrack y,y^{-1} \rbrack $ satisfying $\tau(y) = qy$. We show
the map
\begin{equation}
x^2 \tau - (q+q^{-1})x \tau x + \tau x^2 + (q-q^{-1})^2\tau
\label{eq:andtau}
\end{equation}
vanishes on 
$\K\lbrack y, y^{-1} \rbrack  $. For $f \in  
\K\lbrack y,y^{-1} \rbrack $, the   image of $f$ under   
the map (\ref{eq:andtau}) equals 
\begin{equation}
x^2 - (q+q^{-1})x \tau (x) + \tau (x)^2 + (q-q^{-1})^2
\label{eq:tauaut}
\end{equation}
times
$\tau(f)$.
Observe  the expression (\ref{eq:tauaut}) equals
\beast
(y+y^{-1})^2 - (q+q^{-1})(y+y^{-1})(qy+q^{-1}y^{-1})+ (qy+q^{-1}y^{-1})^2
+(q-q^{-1})^2
\eeast
which equals zero, 
 so the map 
(\ref{eq:andtau}) vanishes on 
$\K\lbrack y,y^{-1} \rbrack $. Similarly   
the map
\begin{equation}
x^2 \tau^{-1} - (q+q^{-1})x \tau^{-1} x + \tau^{-1} x^2 + (q-q^{-1})^2\tau^{-1}
\label{eq:andtauinv}
\end{equation}
vanishes on 
$\K\lbrack y,y^{-1} \rbrack  $. Replacing  $A^*$ by $\D$  in   
(\ref{eq:acom}),
and evaluating the result using (\ref{eq:dddef}) and 
our comments above, we find 
the map (\ref{eq:acom}) is of the form $\omega I$,
where $\omega $ is an element of $\K\lbrack x \rbrack $.
%
We now see $A$ commutes with 
 (\ref{eq:acom}), as desired.
 We have  now shown
(\ref{eq:awTD1}). By
(\ref{eq:awTD1}), 
(\ref{eq:awTD2}), we find $A$ and $A^*$ induce on  
$\K\lbrack x \rbrack $ a module structure for the TD algebra with
parameters 
$\beta = q+q^{-1}$, $\gamma = \gamma^*=0$, $\varrho = -(q-q^{-1})^2$,
$\varrho^*= -abcdq^{-1}(q-q^{-1})^2$. This module is irreducible
in view of Lemma  5.2.

\small
\bibliographystyle{plain}

\bibliography{master}
\normalsize
\medskip
\noindent Paul Terwilliger, Department of Mathematics, University of
Wisconsin, 480 Lincoln Drive, Madison, Wisconsin, 53706, USA \hfil\break
email: terwilli@math.wisc.edu \hfil\break

\end{document}